\documentclass[reqno]{amsart}
\usepackage{enumerate}
\usepackage{amsmath,amssymb,amsthm,amsxtra,amsbsy}
\usepackage{a4}
\usepackage[mathscr]{eucal}

\newtheorem{thm}{Theorem}[section]
\newtheorem{lemma}[thm]{Lemma}
\newtheorem{prop}[thm]{Proposition}
\newtheorem{cor}[thm]{Corollary}
\theoremstyle{definition}

\newtheorem{rem}[thm]{Remark}

\newtheorem{exa}[thm]{Example}

\newcommand{\name}{\textsc}
\renewcommand{\phi}{\varphi}
\renewcommand{\epsilon}{\varepsilon}
\newcommand{\vanish}[1]{\relax}
\newcommand{\sector}[1]{S_{#1}}
\newcommand{\dsector}[1]{\Sigma_{#1}}
\newcommand{\mfrac}[2]{#1/#2}
\newcommand{\strip}[1]{St_{#1}}

\newcommand{\para}[1]{\Pi_{#1}}

\newcommand{\PV}{\mathrm{PV}}
\newcommand{\Mlt}{\mathcal{M}}
\DeclareMathOperator{\Var}{Var}
\newcommand{\res}[1]{ \big|_{#1}}
\newcommand{\DD}{\mathscr{D}}

\newcommand{\RR}{\mathscr{R}}

\newcommand{\resolv}{\varrho}

\DeclareMathOperator{\BIP}{BIP}
\DeclareMathOperator{\Cos}{Cos}

\DeclareMathOperator{\calJ}{\mathcal{J}}
\DeclareMathOperator{\calA}{\mathcal{A}}
\DeclareMathOperator{\calB}{\mathcal{B}}

\DeclareMathOperator{\calE}{\mathcal{E}}
\DeclareMathOperator{\calH}{\mathcal{H}}
\DeclareMathOperator{\calU}{\mathcal{U}}
\DeclareMathOperator{\calX}{\mathcal{X}}
\newcommand{\suchthat}{\,\,|\,\,}
\newcommand{\Sum}[2][\relax]{%
 \ifx#1\relax \sideset{}{_{#2}}\sum 
 \else \sideset{}{^{#1}_{#2}}\sum
 \fi}
\newcommand{\N}{\mathbb{N}}
\newcommand{\Z}{\mathbb{Z}}

\newcommand{\R}{\mathbb{R}}
\newcommand{\C}{\mathbb{C}}
\newcommand{\T}{\mathbb{T}}
\newcommand{\ohne}{\setminus}
\newcommand{\nach}{\circ}
\newcommand{\pfeil}{\longrightarrow}
\newcommand{\tpfeil}{\longmapsto}

\newcommand{\car}{\mathbf{1}}
\DeclareMathOperator{\Hol}{\mathcal{O}}
\DeclareMathOperator{\re}{Re}
\DeclareMathOperator{\im}{Im}
\newcommand{\cls}[1]{\overline{#1}}
\newcommand{\abs}[1]{\left| #1 \right|}

\newcommand{\rand}{\partial}

\newcommand{\twomat}[4]{\left(\begin{array}{cc} #1 & #2 \\
                                                #3 & #4 \end{array}\right)}
\renewcommand{\abs}[1]{\left\vert#1\right\vert}
\DeclareMathOperator{\Lin}{\mathcal{L}}
\newcommand{\norm}[2][\relax]{%
   \ifx#1\relax \ensuremath{\left\Vert#2\right\Vert}
   \else \ensuremath{\left\Vert#2\right\Vert_{#1}}
   \fi}
\makeatletter
\newcommand{\sprod}[2]{\ensuremath{%
  \setbox0=\hbox{\ensuremath{#2}}
  \dimen@\ht0
  \advance\dimen@ by \dp0
  \left(\left.#1\rule[-\dp0]{0pt}{\dimen@}\,\right|#2\hspace{1pt}\right)}}
 \makeatother
\newcommand{\dprod}[2]{\ensuremath{\left<#1,#2\right>}}

\newcommand{\spacefont}{\mathbf}

\newcommand{\BV}{\mathbf{BV}}

\newcommand{\Ct}[2][\relax]{%
 \ifx#1\relax \ensuremath{ {\mathbf{C}}^{\mathbf{#2}}  }
 \else {\mathbf{C}}^{\mathbf{#2}}_{\boldsymbol{#1}}
 \fi}

\newcommand{\Co}[1][\relax]{%
\mathbf{C}_\mathbf{0}^\mathbf{#1}}
\newcommand{\Cb}{\mathbf{C}^\mathbf{b}}

\newcommand{\Bt}[2][\relax]{%
 \ifx#1\relax \ensuremath{ {\mathbf{B}}^{\mathbf{#2}}  }
 \else {\mathbf{B}}^{\mathbf{#2}}_{\boldsymbol{#1}}
 \fi}

\newcommand{\ct}[1][\relax]{%
 \ifx#1\relax \mathbf{c}
 \else \mathbf{c}_{\boldsymbol{#1}}
 \fi}
\newcommand{\loc}{\text{\upshape{\tiny \bfseries loc}}}
\newcommand{\Ell}[2][\relax]{%
   \ifx#1\relax \mathbf{L}^{\boldsymbol{#2}}
   \else \mathbf{L}^{\boldsymbol{#2}}_{\boldsymbol{#1}}
   \fi}
\newcommand{\Wee}[2][\relax]{%
   \ifx#1\relax \mathbf{W}^{\boldsymbol{#2}}
   \else \mathbf{W}^{\boldsymbol{#2}}_{\boldsymbol{#1}}
   \fi}
\newcommand{\Har}[2][\relax]{%
   \ifx#1\relax \boldsymbol{\mathsf{H}}^{\boldsymbol{#2}}
   \else \boldsymbol{\mathsf{H}}^{\boldsymbol{#2}}_{\boldsymbol{#1}}
   \fi}
\newcommand{\eM}{\spacefont{M}}

\newcommand{\Fourier}{\mathcal{F}}
\newcommand{\fourier}[1]{\widehat{#1}}
\newcommand{\tfunk}{\boldsymbol{\mathcal{D}}} 

\newcounter{aufzi}
\newenvironment{aufzi}{\begin{list}{ {\upshape\alph{aufzi})}}{
        \usecounter{aufzi}
        \topsep1ex
        \parsep0cm
        \itemsep1ex
        \leftmargin1cm
        \labelwidth0.5cm
        \labelsep0.3cm
}}
{\end{list}}
\newcounter{aufzii}

\newcounter{aufziii}
\newenvironment{aufziii}{\begin{list}{ {\upshape\arabic{aufziii})}}{
        \usecounter{aufziii}
        \topsep1ex
        \parsep0cm
        \itemsep1ex
        \leftmargin0.8cm
        \labelwidth0.5cm
        \labelsep0.3cm
}}
{\end{list}}
\newcounter{aufziv}

\begin{document}

\title[Transference Principle for General Groups]{A Transference Principle for General Groups and Functional Calculus
on UMD Spaces}
\author{\textsc{Markus Haase}}

\address{Delft Institute of Applied Mathematics, Technical University Delft,
PO Box 5031, 2600 GA Delft, The Netherlands}
\email{m.h.a.haase@tudelf.nl} 

\thanks{Author's e-mail address: {\tt m.h.a.haase@tudelft.nl}}

\date{\today}

\begin{abstract}
{\sffamily Let $-iA$ be the generator of a $C_0$-group 
$(U(s))_{s\in \R}$ on a Banach space $X$, and $\omega > \theta(U)$,
the group type of $U$.
 We prove a transference
principle that allows to estimate $\norm{f(A)}$ in terms of the
$\Ell{p}(\R;X)$-Fourier multiplier norm of $f(\cdot \pm i \omega)$.
If $X$ is a Hilbert space
this yields new proofs of important results of McIntosh and
Boyadzhiev--de Laubenfels. If $X$ is a UMD space, one obtains 
a bounded $H^\infty_1$-calculus of $A$ on horizontal strips. Related
results for sectorial and parabola-type operators follow. Finally 
it is proved  that each generator of a cosine function on a UMD space
has bounded $H^\infty$-calculus on sectors.}


\end{abstract}

\renewcommand{\subjclassname}{\textup{2000} Mathematics Subject Classification}
\maketitle

\section{Introduction}\label{s.intro}

The name ``transference principle'
 was introduced 
by \name{Coifman} and \name{Weiss} in their
influential monograph \cite{CoiWei}, building on earlier
work of \name{Calder\'on} \cite{Cal68}.
The orginal setting is in terms of a locally compact abelian group
acting boundedly on some space $X = \Ell{p}(\Omega,\mu)$, but it was 
observed by \name{Berkson}, \name{Gillespie}
and \name{Muhly} in \cite{BerGilMuh89b}
that the restriction to $\Ell{p}$-spaces was unnecessary and
one could in fact take as $X$ any Banach space. Their result was later on
generalised by \name{Berkson}, \name{Paluzy{\'n}ski} and \name{Weiss}
to so-called {\em transference couples} \cite{BerPalWei99}.
For the group
of real numbers, the ``classical'' transference principle has the following
form (see \cite[Theorem 2.8]{BerGilMuh89b} or 
\cite[Theorem 10.5]{KunWei04} for proofs).

\begin{thm}{\bf (Coifman--Weiss \cite{CoiWei}, 
Berkson--Gillespie--Muhly \cite{BerGilMuh89b})}\label{tra.t.bdg}\\
Let $U$ be a $C_0$-group on a Banach space $X$,
such that $\norm{U(t)} \le M$ for all $t \in \R$. Then
\[ \norm{ \int_\R U(s)x\, \mu(ds)} \le M^2 \norm{L_\mu}_{\Lin(\Ell{p}(\R;X))}
\norm{x}_X 
\]
for all $x\in X, \mu \in \eM(\R)$, $1\le p < \infty$.
Here, $L_\mu$ denotes the convolution operator $L_\mu := 
(f \tpfeil \mu \ast f)$.
\end{thm}

So  ``transference''
means that certain averages over the representation of the group can be
estimated by the norm of an associated convolution operator. In this
manner one can for example prove that restrictions to $\Z$ of Fourier
multipliers on $\R$ are Fourier multipliers on $\T$, but there are plenty of
other examples and applications of transference, see 
\cite{CoiWei77, BerGil87, BerGil89, BerGil90, SorWeiZal90, AsmBerGil92b, 
BerGil94, AsmKel96, Bel99}. 

Although the transference result in itself is valid in any Banach
space, to be applicable one usually needs special conditions. The by far most
useful of these 
seems to be the so-called UMD-property (see Section \ref{s.trans} below).  
It were \name{Cl\'ement} and \name{Pr\"uss}
 who observed in \cite{ClePru90} the usefulness of the vector-valued
(in particular: UMD-valued) transference principle
for treating problems in evolution equations. The link is via the
group of imaginary powers of a sectorial operator, and is based on
the 1987 paper \cite{DorVen87} of \name{Dore} and \name{Venni} 
which established the relation between bounded imaginary powers and
the maximal regularity property of a sectorial operator. A few years later,
when \name{McIntosh}'s notion of a bounded $H^\infty$-calculus 
(from \cite{McI86}) had become adopted by the evolution equations community,
\name{Hieber} and  \name{Pr\"uss} 
used  Theorem \ref{tra.t.bdg}
to derive the following result.

\begin{thm}{\bf (Hieber--Pr\"uss \cite{HiePru98})}\label{int.t.hiepru}\\
Let $-iA$ generate a bounded $C_0$-group on a UMD space $X$. Then
for every $\phi \in (0, \pi/2)$ the operator $A$ has a bounded
$H^\infty(\dsector{\phi})$-calculus.
\end{thm}

\noindent
(The symbol $\dsector{\phi}$ denotes the double sector of angle $\phi$, see
 below.) For proofs see also \cite[Theorem 10.7]{KunWei04} and
Remark \ref{tra.r.hiepru} below. 

Around the year 2000 the notion of $R$-boundedness was fully developed
and its importance for  operator-valued Fourier
multiplier theorems and the maximal regularity problem had been
recognised by many authors; key steps were
the paper \cite{CledePSukWit00} and, of course, \name{Weis}' paper
\cite{Wei01}. In \cite{ClePru01}, using the notion of $R$-boundedness, 
\name{Cl\'ement} and \name{Pr\"uss}
extended 
the transference principle and Theorem \ref{int.t.hiepru} 
to the operator-valued setting  and gave applications to
the maximal regularity problem.
See \cite{Wei01, ClePru01, DenHiePru03a, KunWei04} for more information.

Despite the wide range of applications, it is clear that confining oneself
to bounded representations of groups is a major restriction; however, it
seems that no transference result for unbounded groups is available yet
in the literature. In this article we shall fill this gap (Theorem 
\ref{tra.t.tra2}). Moreover, we shall apply the transference principle
to obtain a generalisation of the Hieber--Pr\"uss theorem to
unbounded groups (Theorem \ref{ex.t.ven}), and to obtain
boundedness results for  certain functional calculi of
unbounded operators (Theorem \ref{cos.t.mci}); these results  are
the UMD-analogues of the Hilbert space theorems of \name{McIntosh} \cite{McI86}
and of \name{Boyadzhiev} and \name{de Laubenfels} \cite{BoydeL94}.
Other applications are to generators of cosine functions on UMD spaces
(Theorem \ref{cos.t.fcpar}, Corollary \ref{cos.c.fcsec}).
Our results can also be used to obtain
new proofs for important results of \name{Fattorini} 
\cite[Theorem 3.16.7]{ABHN}, \name{Dore} and \name{Venni} \cite{DorVen87} 
and \name{Monniaux} \cite{Mon99}, see \cite{Haa07b}. The
transference method was also used successfully in \cite{Haa08b}
to resolve an old problem concerning bounded cosine functions  on UMD spaces.

\smallskip

The paper is organised as follows. In Section \ref{s.fc} we provide
without proofs 
the necessary material on functional calculus.
In Section \ref{s.trans} we prove the main transference result
(Theorem \ref{tra.t.tra2}). Then we introduce UMD spaces and use
\name{Bourgain}'s UMD-valued multiplier result of Mikhlin type to
derive our main result on functional calculus
(Theorem \ref{tra.t.fc}). In Section \ref{s.ex} we discuss some
examples, in  particular, we describe a class of functions $f$
such that the operator $f(A)$ is explicitly given by a
principal-value integral (Theorem \ref{ex.t.bv}). 
We also give the generalisation of the
Hieber--Pr\"uss result to unbounded groups. Section \ref{s.cos}
is devoted to the consequences of our main results in the context
of sectorial and parabola-type operators.

\bigskip

\noindent
{\em Notation}\\
We usually consider (unbounded) closed
operators $A,B$ on a Banach space $X$.  By $\Lin(X)$ we denote
 the set of all {\em bounded} (fully-defined) operators on $X$.
The domain and the range
of a general  operator $A$ are denoted by $\DD(A)$ and $\RR(A)$, respectively. 
Its resolvent is $R(\lambda,A) = (\lambda-A)^{-1}$, and $\resolv(A)$
denotes the set of $\lambda \in \C$ where $R(\lambda,A) \in \Lin(X)$.
Its complement $\sigma(A) = \C\ohne \resolv(A)$ is the {\em spectrum}.
For given $\omega > 0$ we define
\[ \strip{\omega} := \{ z\in \C \suchthat \abs{\im z } < \omega\}
\]
to be the horizontal strip of width $2\omega$. If $\omega = 0$
we define $\strip{0} := \R$. Accordingly
\[ \para{\omega}
:= \{ z^2 \suchthat z\in \strip{\omega} \} 
\]
is the horizontal parabola, symmetric about $\R$ and unbounded to the right.
Furthermore,
for $\omega \in [0, \pi]$  
\[ \sector{\omega} := \{ e^z  \suchthat z\in \strip{\omega} \}
= \begin{cases}
\{ z \not= 0 \suchthat \abs{\arg z} < \omega\}, & \omega \in (0, \pi]\\
    (0, \infty), & \omega = 0
\end{cases}
\]
is the horizontal sector of angle $2\omega$, symmetric about the
positive real axis. For $\omega \in [0, \pi/2)$
\[ \dsector{\omega} := \sector{\omega} \cup -\sector{\omega}
\]
is the horizontal double sector.

For each open subset $\Omega \subset \C$
we denote by $H^\infty(\Omega)$
the Banach algebra of bounded holomorphic functions on $\Omega$.
If $\Omega$ is an arbitrary locally compact space, then 
the set of complex regular Borel measures on $\Omega$ is 
denoted by $\eM(\Omega)$.
The {\em Fourier transform} of a tempered distribution $\Phi$ on $\R$ is 
denoted by $\Fourier(\Phi)$ or $\fourier{\Phi}$. We often write
$s$ and $t$ (in the Fourier image) to denote the real coordinate, e.g.
$\sin t/ t$ denotes the {\em function} $t \mapsto  \sin t/ t$. 

Let $X$ be a Banach space. For a finite measure $\nu \in \eM(\R)$ we
denote by
\[ L_\nu := (f \to f \ast \nu): \Ell{p}(\R;X) \to \Ell{p}(\R;X)
\]
the convolution operator on the $X$-valued $\Ell{p}$-space.

\section{Functional Calculus Prelims}\label{s.fc}

In this section we provide necessary facts on functional calculus. Most of the 
results are well-known, cf. \cite{PruSoh90}, \cite{HiePru98}
and the comments in \cite{HaaFC}.

\medskip

Let $-iA$ be the generator of a $C_0$-group $(U(s))_{s\in \R}$ on a Banach
space $X$. It is an elementary fact that there exist constants $M\ge 1, 
\omega \ge 0$ such that 
\[ \|U(s)\| \le M e^{\omega \abs{s}}, \quad \quad (s\in \R).
\]
The infimum of all such $\omega \ge 0$ is called the {\em group type}
of $U$ and is (here) denoted by $\theta(U)$.
Abstract semigroup theory yields that the resolvent of $A$ satisfies the
estimate
\[
\norm{R(\lambda,A)} \le \frac{M}{\abs{\im \lambda} - \omega} \quad \quad
(\abs{\im \lambda} > \omega),
\]
hence $A$ is a so-called strong strip-type
operator of type $\omega_{sst}(A) \le \omega$, as defined 
in \cite[Section 4.1]{HaaFC}.
There is a natural holomorphic functional calculus associated with
such operators: as a first step one uses the Cauchy formula to 
define
\[ f(A) :=  \frac{1}{2\pi i} \int_\Gamma f(z) R(z,A)\, dz,
\]
where $f$ is a holomorphic function on a strip in the class
$\calE(\theta)$ defined by
\[ \calE(\theta) := \left\{ f \in H^\infty(\strip{\theta}) \suchthat
f(z) = O\big(\abs{z}^{-2}\big)\,\, (\re z \to \pm \infty)\right\},
\] 
for some $\theta > \omega$.
The contour $\Gamma$ is the positively oriented
boundary of a smaller strip $\strip{\omega'}$, 
$\omega' \in (\omega, \theta)$ being arbitrary. This yields an algebra 
homomorphism of the algebra $\calE(\theta)$  into the
algebra of bounded operators on $X$. In a second step, by so-called
regularisation, one defines $f(A)$ for a much wider class of functions:
\[ f(A) := e(A)^{-1} (ef)(A),
\]
where $e\in \calE(\theta)$ is such that also $ef \in \calE(\theta)$ and
$e(A)$ is injective. The function $e$ is called a {\em regulariser} for $f$,
and the definition of $f(A)$ is independent of the chosen regulariser.
For example, if $f\in H^\infty(\strip{\theta})$, then $f$ is regularisable by
any function $e(z) = (\lambda - z)^{-2}$, where $\abs{\im \lambda} > 
\theta$. 
Details of the construction as well as a listing of all the formal 
properties of the so constructed functional calculus can be found in 
\cite[Chapter 1 and 4]{HaaFC}.

Of particular importance in the theory of functional calculus is the so-called
convergence lemma. It goes back to \name{McIntosh} \cite{McI86} in  the 
sectorial case; for a proof see  \cite{Haa07b} or
\cite[Prop. 5.1.4]{HaaFC}.

\begin{prop}{\bf (Convergence Lemma)}\label{fc.p.cl}\\
Let $A$ be a strong strip-type operator on the Banach space $X$, 
with dense domain. Let  $\theta> \omega_{sst}(A)$, and let
$(f_\iota)_{\iota\in \calJ}$ be a net of holomorphic functions on the
strip $\strip{\theta}$, satisfying
\begin{aufziii}
\item $\sup \{ \abs{f_\iota(z)} \suchthat z\in \strip{\theta}, 
\iota\in \calJ\} < \infty$;
\item $f_\iota(z) \to f(z)$ for every $z\in \strip{\theta}$;
\item $\sup_\iota \norm{f_\iota(A)} < \infty$.
\end{aufziii}
Then $f(A) \in \Lin(X)$ and $f_\iota(A) \to f(A)$ strongly.
\end{prop}

In the case that $-iA$ generates a $C_0$-group there
is a convenient tool to 
identify functions $f$ such that $f(A)$ is a bounded operator. First of all,
consider for $s\in \R$ the function 
$e^{-isz}$, which is bounded on every horizontal strip.
We clearly expect $e^{-isz}(A) = U(s)$ for all $s\in \R$.  More generally,
let $\omega \ge 0$ be fixed such that $\norm{U(s)} \le M e^{\omega\abs{s}}$
for some $M\ge 1$ and all $s\in \R$, and 
let $\mu$ be a (complex) Borel measure on $\R$ satisfying
\begin{equation}\label{fc.e.mo}
 \norm{\mu}_{\eM_\omega} := \int_\R e^{\omega\abs{s}} \, \abs{\mu}(ds) < 
\infty.
\end{equation}
Then one can set
\[ T_\mu x := \int_\R U(s)x\, \mu(ds)\quad \quad (x\in X).
\]
Clearly, the set
$\eM_\omega(\R) := \{ \mu \suchthat \text{(\ref{fc.e.mo}) holds}\}$
is a Banach algebra with respect to convolution, and the map
$(\mu \tpfeil T_\mu): \eM_\omega(\R) \pfeil \Lin(X)$ is
a homomorphism of algebras, called the
{\em Phillips calculus}.
Of course  we expect $T_\mu = f(A)$ 
where $f$ is the {\em Fourier-Stieltjes}
transform
\[ 
f(z) = \fourier{\mu}(z) 
:= \int_\R e^{-isz} \, \mu(ds)\quad \quad (z\in \strip{\theta})
\]
of $\mu$. Note that $\fourier{\mu} \in 
H^\infty(\strip{\omega}) \cap \Cb(\cls{\strip{\omega}})$.
Here is the precise result, a proof of which is in \cite{Haa07b},
see also \cite[Section 2]{HiePru98}.

\begin{lemma}\label{fc.l.pc}
Let $X$, $A$, and $U$ be as above, and let $\theta > \omega$.
\begin{aufzi}
\item Each function
$f \in \calE(\theta)$ arises as a Fourier--Stieltjes transform, namely
\[ f = \fourier{g} \quad \text{with}\quad 
g(s) := \frac{1}{2\pi} \int_\R f(t) e^{ist}\, dt
\quad \quad (s\in \R).
\]
One has $g\in \Co{}(\R)$ and $\int \abs{g(s)} 
e^{\alpha \abs{s}} \, ds < \infty$ for all $\alpha \in [0, \theta)$; in 
particular, one has  $g(s)ds \in \eM_\omega(\R)$.
\item Let $\mu \in \eM_\omega(\R)$, and suppose that 
  $f:= \fourier{\mu}$ extends to a holomorphic function on $\strip{\theta}$
such that $f(A)$ is defined.  Then
$f(A) = T_\mu \in \Lin(X)$ and 
\[ \sup_{t \in \R} \norm{f(t+A)} \le M \norm{\mu}_{\eM_\omega}.
\] 
\end{aufzi} 
\end{lemma}

The formulation of part b) in the above proposition is due to the
fact that in our orginal set-up of the functional calculus, functions
had to be defined on strips strictly larger than the spectral strip.
By Lemma \ref{fc.l.pc} we can extend the orginal definition and write
\[ f(A)x := T_\mu x := \int_\R U(t)x\, \mu(dt) \quad \quad (x\in X)
\]
whenever $f = \fourier{\mu}, \mu \in \eM_\omega(\R)$. This will induce
also a compatible extension of the unbounded functional calculus, 
see \cite[Proposition 1.2.7]{HaaFC}.

However, as far as {\em bounded} operators are concerned,
one cannot go beyond Fourier transforms
of measures $\mu \in \eM_\omega(\R)$ in general. Indeed, 
if $f(A)$ is bounded in the special case of  $-iA= d/dt$ generating
the shift group on  
\[
X := \Ell[\omega]{1}(\R) := \Big\{ f \in \Ell[\loc]{1}(\R) \suchthat
\int_\R f(t) e^{\omega \abs{t}}\, dt < \infty\Big\},
\]
then actually $f = \fourier{\mu}$ for some $\mu \in \eM_\omega(\R)$
(see for example  \cite[Proposition 2.3]{Haa07b}). 
One can put this remark in the form of a transference principle.

\begin{prop}{\bf (Transference Principle, $\Ell{1}$-version)}\\
Let $-iA$ generate $C_0$-group $U$ on  a Banach space $X$, and let
$M\ge 1, \omega \ge 0$ be such that $e^{-\omega \abs{t}} \norm{U(t)} \le M$
for all $t\in \R$. Then
\[ \norm{f(A)}_{\Lin(X)} \le 
M \norm{f(id/dt)}_{\Lin(\Ell[\omega]{1}(\R))} 
\]
for all $f = \fourier{\mu}$, $\mu \in \eM_{\omega}(\R)$.
\end{prop}

\begin{proof}
Let $\mu \in  \eM_{\omega}(\R)$, define $f := \fourier{\mu}$.
If $A_0 := id/dt$ on $\Ell[\omega]{1}(\R)$, it is easily seen
that (1) $f(A_0)$ is convolution with the measure $\mu^\sim$,
defined by $\mu^\sim(B) = \mu(-B)$, and that
(2) the norm of this convolution operator is exactly 
$\norm{\mu}_{\eM_\omega}$ (see \cite[(2.2)]{Haa07b}).
But it is trivial that
\[ \norm{f(A)} = \norm{T_\mu} \le M \norm{\mu}_{\eM_\omega}
\]
and so we are done.
\end{proof}

Let us mention that although this first example of a transference
principle is fairly elementary, it has important consequences.
See e.g.~\cite{BaeKov06pre} for the case of bounded groups.

It is known for a long time that if
$-iA$ generates a group on a UMD space $X$,
$f(A)$ is bounded even for  functions $f$ arising
from certain principal value distributions (see also Section
\ref{s.ex} below.) In fact, this is the core of the results of 
Dore--Venni, Monniaux and Fattorini, but this  
was fully brought to light only recently in \cite{Haa07b}. 
(We are oversimplifying
here, but a more detailed discussion would lead us too far astray.)
In the following we shall extend the results from \cite{Haa07b} 
towards a full bounded functional calculus.

\section{The Transference Principle}\label{s.trans}

In the effort to obtain an  analogue for unbounded groups of 
the classical Coifman--Weiss transference principle 
(Theorem \ref{tra.t.bdg}) a first 
progress was made recently in \cite[Theorem 3.1]{Haa07b}.
We  reproduce the result here, because its proof is short and instructive.

\begin{thm}{\bf (Transference Principle, fixed compact support 
\cite{Haa07b})}%
\label{tra.t.tra1}\\
Let $p \in [1, \infty)$, and let $U$ be a $C_0$-group on a Banach space $X$ 
Define $M := \sup_{s\in [-2,2]} \norm{U(s)}$. 
Then 
\[
\norm{ \int_{[-1,1]} U(s)x\, \mu(ds) } \le 2^{1/p} M^2  
\norm{L_\mu}_{\Lin(\Ell{p}(\R;X))} \norm{x} \quad \quad (x\in X)
\]
for all $\mu \in \eM[-1,1]$. 
\end{thm}

\begin{proof}
Let $\abs{t} \le 1$. Then we write 
\[
T_\mu x= \int_{-1}^1 U(s)x\, \mu(ds) = U(t) \int_{-1}^1 U(s-t)x\, \mu(ds)
= U(t)(f \ast \mu)(t)
\] 
where $f(s) = \car_{[-2,2]}(s) \,U(-s)x$.
Hence
\[ \norm{T_\mu x} = \norm{\frac{1}{2} \int_{-1}^1 U(t) (f\ast \mu)(t)\, dt}
\le 2^{\frac{1}{p'}-1}M \norm{f \ast \mu}_{\Ell{p}(R;X)}
\le 2^{\frac{1}{p}} M^2 \norm{L_\mu} \norm{x},
\] 
where $\norm{L_\mu} = \norm{L_\mu}_{\Lin(\Ell{p}(\R;X))}$.
\end{proof}

In Theorem \ref{tra.t.tra1} 
the actual growth type of the group is irrelevant since
only measures $\mu$  with a fixed compact support are considered. If we allow
unbounded support, we have to modify the proof. However, 
the strategy will be the same: we first convolve the 
(slightly modified) measure $\mu$ with some function in $\Ell{p}$,
then we integrate (against a ``test function'') and in the end
obtain the operator $T_\mu$. 
Given $\mu \in \eM_\omega(\R)$ we define the finite measure
$\mu_\omega$ by
\[ \mu_\omega(ds) = \cosh(\omega s)\mu(ds).
\]
Then the result is as follows.

\begin{thm}{\bf (Transference Principle)}%
\label{tra.t.tra2}\\
Let $0 \le \omega_0 <  \omega$, and 
$p\in [1, \infty)$. Then there is a constant 
$C = C(p, \omega_0, \omega)$ such that the following holds:
If $U$ is a $C_0$-group on a Banach space $X$ such that
$\norm{U(s)} \le M \cosh(\omega_0 s)$ for all $s\in \R$ and some
$M \ge 1$, then  
\[
\norm{ \int_\R U(s)x\, \mu(ds) } \le C M^2  
\norm{L_{\mu_\omega}}_{\Lin(\Ell{p}(\R;X))} \norm{x} \quad \quad (x\in X)
\]
for all $\mu \in \eM_\omega(\R)$.
\end{thm}

\begin{proof}
Fix $\alpha > \omega$ and $x\in X$, and consider the function
$f$, defined by
\[ f(s) := [\cosh(\alpha s)]^{-1}U(-s)x, \quad \quad (s\in \R).
\] 
Clearly $f \in \Ell{p}(\R;X)$ and 
\[ \norm{f}_{\Ell{p}(\R;X)}
\le M c_1 \norm{x}
\]
with $c_1 := \norm{\cosh(\omega_0 s)/ \cosh(\alpha s)}_{\Ell{p}(\R)}$.
Hence
\[ \norm{f \ast \mu_\omega}_{\Ell{p}} = \norm{L_{\mu_\omega} f}_{\Ell{p}}
\le M \norm{L_{\mu_\omega}}_{\Lin(\Ell{p}(\R;X))} c_1 \norm{x}.
\]
Note that for every $t\in \R$
\begin{align*}
 (f \ast \mu_\omega)(t) & = \int_\R  U(s-t)x 
\frac{\cosh(\omega s)}{\cosh(\alpha(s-t))}\, \mu(ds) \\
& =  U(-t) \int_\R  U(s)x 
\frac{\cosh(\omega s)}{\cosh(\alpha(s-t))}\, \mu(ds).
\end{align*}
Let $\phi$ be a scalar function such that 
$\phi(s)= O( \cosh(\omega s)^{-1})$ as 
$\abs{s} \to \infty$. Then
\begin{align*}
T & := \int_\R \phi(t) U(t)[f \ast \mu_\omega](t)\, dt
\\ & = \int_\R  [\phi \ast \cosh(\alpha \,\cdot)^{-1}](s)\, \cosh(\omega s) \,U(s)x\, \mu(ds)
\end{align*}
by Fubini's theorem. If we can choose $\phi$ such that
\begin{equation}\label{tra.e.dec}
\left[\phi \ast   \cosh(\alpha\, \cdot)^{-1}\right](s)
= [\cosh(\omega \,s )]^{-1}\quad \quad (s\in \R)
\end{equation}
 then
$T = T_\mu$ and
\[ \norm{T_\mu} = \norm{T} \le M 
\norm{ \phi \cosh(\omega_0\, \cdot)}_{\Ell{p'}}
M \norm{L_{\mu_\omega}}_{\Lin(\Ell{p}(\R;X))} c_1 \norm{x}.
\]
To determine such a $\phi$ we take Fourier transforms on both sides
of (\ref{tra.e.dec}). It is known that the function
$\cosh^{-1}$ is almost its own
Fourier transform. More precisely, one has
\[ \int_\R \frac{e^{-isz}}{\cosh(\omega s)}\, ds =
\frac{\pi / \omega}{\cosh((\pi/2\omega) z)} \quad \quad 
(\abs{\im z} < \omega)
\]
(see \cite[p.81]{SteSha2} for a proof).
Hence we look for a function $\phi$ that satisfies
\[ \fourier{\phi}(z)\,\, \frac{\pi / \alpha}{\cosh((\pi/2\alpha) z)} =
\frac{\pi / \omega}{\cosh((\pi/2\omega) z)} 
\]
that is
\[ \phi(t) = \frac{\alpha}{\omega}\,
\Fourier^{-1}\left(  
\frac{\cosh((\pi/2\alpha) \,\cdot)  }{\cosh((\pi/2\omega) \,\cdot)} 
\right)(t) =
\frac{2\alpha}{\pi} \cos\left(\frac{\pi \omega}{2\alpha}\right)
\frac{\cosh(\omega t)}{\cos(\pi \omega/ \alpha) + \cosh(2\omega t)}
\]
(computed from the second formula of \cite[p.36]{Obe}). For example,
taking $\alpha = 2\omega$ the function
\[ \phi(t) = \frac{\sqrt{8}\omega}{\pi}\, 
\frac{\cosh(\omega t)}{\cosh(2\omega t)}\quad \quad(t\in \R)
\]
will do.
\end{proof}

\begin{rem}\label{tra.r.tra2}
It is worthwhile to compare
Theorem \ref{tra.t.tra2} with Theorem \ref{tra.t.bdg}. In the proof of the
former it is essential that $\omega$ 
is strictly larger than $\omega_0$,
the exponential growth type of the group. This matches the experience
that bounded groups behave much better than general ones. 
We do not expect the conclusion of Theorem \ref{tra.t.tra2} to hold
for $\omega = \omega_0$.
\end{rem}

As in the classical case, one can phrase the transference principle
in terms of norms of Fourier multipliers. Namely, one has
\[ 2 \mu_\omega(ds) = 2 \cosh(\omega s)\mu(ds) = 
e^{\omega s}\mu(ds) + e^{-\omega s}\mu(ds),
\] 
and since
\[ \Fourier( e^{\omega \,\cdot} \mu) = \fourier{\mu}(\,\cdot - i\omega)
\quad \text{and} \quad 
\Fourier( e^{-\omega \,\cdot} \mu) = \fourier{\mu}(\,\cdot + i\omega),
\]
the convolution operator $L_{\mu_\omega}$ from above is the average
of the two Fourier multiplier operators with symbols
$\fourier{\mu}(\,\cdot \pm i\omega)$. Denoting the space of
bounded $\Ell{p}(\R;X)$-Fourier multipliers with $\Mlt_p(X)$, we
obtain the following.

\begin{cor}\label{tra.c.fm}
Let $0 \le \omega_0 <  \omega$  and 
$p\in [1, \infty)$. Then there is a constant 
$C = C(p, \omega_0, \omega)$ such that the following holds:
If $U$ is a $C_0$-group on a Banach space $X$ such that
$\norm{U(s)} \le M \cosh(\omega_0 s)$ for all $s\in \R$ and some
$M \ge 1$, then   
\[
\norm{ \int_\R U(s)x\, \mu(ds) } \le C M^2  
\Big( \norm{\fourier{\mu}(\,\cdot + i\omega)}_{\Mlt_p(X)} 
+ \norm{\fourier{\mu}(\,\cdot - i\omega)}_{\Mlt_p(X)} \Big) 
\norm{x}
\]
for all $x\in X$ and all $\mu \in \eM_\omega(\R)$.
\end{cor}

As in the classical case, the transference principle shows its full strength
in a situation when one actually knows something about 
$\Ell{p}(\R;X)$-Fourier multipliers. This is the case when the Banach
space has the so-called UMD property.

\smallskip

A Banach space $X$ has the UMD property if 
every $X$-valued $\Ell{p}$-martingale has unconditional differences.
This notion was introduced by \name{Burkholder} \cite{Bur81a} and it 
is independent of $p \in (1, \infty)$. \name{Burkholder} \cite{Bur83} and 
\name{Bourgain} \cite{Bou83a} showed 
that the UMD property is in fact equivalent to the boundedness
of the Hilbert transform on $\Ell{2}(\R;X)$. 
To make this explicit, consider
the truncated Hilbert transform
\[
\calH_\epsilon f(s) := \int_{\abs{t}\ge \epsilon} f(s-t)\frac{dt}{t} \quad \quad
(f\in \Ell{2}(\R;X)).
\]
Then $X$ is UMD if and only if $(T_\epsilon)_\epsilon$ is uniformly
bounded in $\Lin(\Ell{2}(\R;X))$ if and only if
$(\calH_\epsilon)_\epsilon$ converges as $\epsilon \searrow 0$ 
strongly in $\Lin(\Ell{2}(\R;X))$
to a bounded operator $\calH$. The operator $\calH$  is called the
Hilbert transform. On UMD spaces we have the following Mikhlin-type
multiplier theorem \cite{Zim89}.

\begin{lemma}{\bf (Mikhlin, UMD-valued \cite{Zim89})}\label{tra.l.mik}\\
Let $X$ be a UMD space and $p \in (1, \infty)$. Then there is a constant $C_p$
such that every 
$m \in \Ct{1}(\R\ohne\{0\})$ with 
\[ c_m := \sup_{t\in \R\ohne\{0\}} \abs{m(t)} + 
\sup_{t\in \R\ohne\{0\}} \abs{tm '(t)}  <\infty.
\]
is a bounded $\Ell{p}(\R;X)$-multiplier such that
$\norm{m}_{\Mlt_p(X)} \le C_p c_m$.
\end{lemma}

Using this result, we can now state and prove the main result on 
functional calculus. For $\theta > 0$ define
\[ H^\infty_1(\strip{\theta}) := \{ f \in \Hol(\strip{\theta}) \suchthat
f, z f' \in H^\infty(\strip{\theta})\}.
\]
This set is a Banach algebra with the norm
\[
\norm{f}_{H^\infty_1} := \sup_{z\in \strip{\theta}} \abs{f(z)} + 
\abs{z f'(z)}. 
\]
Every elementary rational function $(\lambda - z)^{-1}$, 
$\abs{\im \lambda} > \theta$, belongs to $H^\infty_1(\strip{\theta})$.

\begin{thm}{\bf (Functional Calculus)}\label{tra.t.fc}\\
Let $X$ be a UMD space, and let $-iA$ be the generator
of a strongly continuous group $U= (U(s))_{s\in \R}$ on $X$. 
Let $\theta > \theta(U)$. 
Then $f(A+r)$ is bounded for every $f\in H^\infty_1(\strip{\theta})$ and
every $r\in \R$, and 
there is a constant $c > 0$ such that
\[
\norm{f(A+ r)} \le c\,  \big(\norm{f}_\infty + \norm{z f'(z)}_\infty\big) \quad \quad
(f\in H^\infty_1(\strip{\theta}), r\in \R).
\]
\end{thm}

\begin{proof}
Let $M \ge 1$ and $\omega_0\in [0,\theta)$ such that $\norm{U(s)} \le M 
\cosh(\omega_0 s)$, $s\in \R$, and choose $\omega \in (\omega_0, \theta)$.
In a first step, take $f\in H^\infty_1(\strip{\theta}) \cap 
\calE(\strip{\theta})$. Then by Lemma \ref{fc.l.pc} there is 
a function $g$ on $\R$ such that $g(s)ds \in \eM(\R)$ and $\fourier{g} = f$.
Fixing $x\in X$ and applying Corollary \ref{tra.c.fm} yields
\begin{align*}
 \norm{f(A+r)x} & = \norm{ \int_\R g(s) e^{-isr}U(s)x\, ds}
\\ & \le C M^2 \norm{x}\, \left( \norm{f(\,\cdot + i\omega)}_{\Mlt_2(X)} 
+ \norm{f(\,\cdot - i\omega)}_{\Mlt_2(X)} \right)
\end{align*}
which by Lemma \ref{tra.l.mik} can be estimated by
\begin{align*}
\dots & \le C' M^2  \norm{x} \left( \norm{f}_\infty + 
\sup_{t \in \R} \abs{t f'(t + i \omega)} + \abs{t f'(t - i \omega)}
\right) 
\\ & \le C' M^2 \norm{x} \norm{f}_{H^\infty_1(\strip{\theta})}.
\end{align*}
To complete the
proof of the theorem, we employ the Convergence Lemma 
(Proposition \ref{fc.p.cl}).
Let $\tau_n(z) := in(in - z)^{-1}$ for $n \in \N$ large.
Then $\tau_n \in H^\infty_1$,
$\sup_n \norm{\tau_n}_{H^\infty_1} < \infty$ and $\tau_n \to 1$ uniformly on 
compacts. 
Let $f\in H^\infty_1$ be arbitrary and define $f_n := f \tau_n^2$.
Then $k := \sup_n \norm{f \tau_n}_{H^\infty_1}
< \infty$, in particular: $\sup_n \norm{f\tau_n}_\infty < \infty$,
and $f \tau_n \to f$ uniformly on compacts. Clearly $f_n \in 
\calE(\strip{\theta})$, so we know already that
\[ \norm{f_n(A+r)} \le c' \norm{f_n}_{H^\infty_1(\strip{\theta})} \le c' k
\]
independent of $n \in \N$ and $r\in \R$. Applying the  Convergence Lemma 
yields that
$f(A+r)\in \Lin(X)$ and 
\[
\norm{f(A+r)}\le C' M^2\, \big(  \limsup_n 
\norm{\tau_n^2}_{H^\infty_1}\big)\,  \norm{f}_{H^\infty_1}
\]
for all $r\in \R$. 
\end{proof}

If $X$ happens to be a Hilbert space, one obtains a much better result. 
In fact, instead of Lemma \ref{tra.l.mik} one can use Plancherel's theorem
and estimate 
\[ \norm{L_{\mu_\omega}} \le 
\frac{1}{2} \norm{ \fourier{\mu}( \,\cdot + i \omega)}_\infty
+ \frac{1}{2} \norm{ \fourier{\mu}( \,\cdot - i \omega)}_\infty 
\le \norm{\fourier{\mu}}_{H^\infty}.
\]
By using the Convergence Lemma as
in the previous proof,  this  leads to the following.

\begin{cor}{\bf (Boyadzhiev--de Laubenfels \cite{BoydeL94})}\label{tra.c.bdl}\\
Let $-iA$ be the generator of a $C_0$-group $U$ on a Hilbert space $H$. Then
for every $\theta > \theta(U)$, the natural $H^\infty(\strip{\theta})$-calculus
is bounded. 
\end{cor}

This theorem was originally proved in \cite{BoydeL94}, but subsequently
reproved in \cite{Haa03a} and \cite{Haa04a}, 
cf. also \cite[Section 7.2]{HaaFC}.

\begin{rem}\label{tra.r.hiepru}
The proof of Theorem \ref{tra.t.fc} 
carries over to a proof of the 
Hieber--Pr\"uss Theorem \ref{int.t.hiepru}; one has to use
Theorem \ref{tra.t.bdg} instead of Theorem \ref{tra.t.tra2}, and
Lemma \ref{ex.l.embed} below.
\end{rem}

\section{Some Classes of Examples}\label{s.ex}

We are going to discuss some classes of functions $f\in H_1^\infty$.
The first avoids involving derivatives.

\begin{lemma}
Let $\omega > 0$ and let $f \in H^\infty(\strip{\omega})$. 
If there exists $a ,b \in \C$ such that
$f(z) -a = O(z^{-1})$ as $\re\ z \to +\infty$ and 
$f(z) -b = O(z^{-1})$ as $\re\ z \to -\infty$, then 
$f \in H^\infty_1(\strip{\theta})$ for every $\theta \in (0, \omega)$.
\end{lemma}

\begin{proof}
Define $g(z) = zf(z) - az$. This is bounded on the half-strip
$(\re z \ge 0, \abs{\im z} < \omega)$. Cauchy's formula yields
\[ g'(z) = \frac{1}{2\pi i} \int_\Gamma \frac{g(w)}{(w-z)^2}\, dw,
\]
where $\Gamma$ is the positively oriented boundary of
a right half-strip within the orginal one and $z$ is within this smaller
half-strip. Consequently, 
$g'$ is uniformly bounded on every even smaller half-strip. But
$g'(z) = f(z) + zf'(z) -a$, and so $zf'(z)$ is bounded on that half-strip.
Analogously, $zf'(z)$ is bounded on left half-strips, and so we 
conclude that $z f'(z)$ is bounded on whole strips 
$\strip{\theta}$, $\theta \in (0,\omega)$.
\end{proof}

We now turn to a class of examples
where one actually has a representation of $f(A)$ as a principal value
integral. 
For an even function
$g \in \Ell{1}(-1,1)$  (i.e., $g(t) = g(-t)$) we define 
the distribution $\PV-g(s)/s$ by the formula
\[
\dprod{\PV-\frac{g(s)}{s}}{\phi} :=
\lim_{\epsilon \searrow 0} \int_{\epsilon < \abs{s} < 1} 
g(s)\phi(s)\frac{ds}{s} = \int_0^1 g(s) \frac{\phi(s) - \phi(-s)}{s} \, ds
\]
for $\phi \in \tfunk(\R)$. Then it is clear that
\[ \abs{\dprod{\PV-\frac{g(s)}{s}}{\phi}} \le \norm{g}_{\Ell{1}(-1,1)}
\norm{\phi'}_\infty
\]
whence $\PV-g(t)/t$ is in fact a distribution of first order.
The proof of the following lemma is easy (see \cite{Haa07b}).

\begin{lemma}\label{ex.l.pvd}
Let $g \in \Ell{1}(-1,1)$ be even and define $h := PV-g(t)/t$ as above.
Then the following assertions hold.
\begin{aufzi}
\item  $h$ is an odd distribution.
\item  Its Fourier transform is
\[ \fourier{h}(z) = \PV-\int_{-1}^1 g(t)e^{-isz} \, \frac{ds}{s}
= (-2i) \int_0^1 \frac{\sin(sz)}{s} g(s) \, ds \quad \quad (z\in \C)
\]
\item One has 
$\displaystyle \frac{d}{dz} \fourier{h}(z) = (-i) \fourier{g}(z)$,
$z\in \C$, and $\fourier{h}(0) = 0$.
\end{aufzi}
\end{lemma}

\vanish{
\begin{exa}
Take $g(s) = 1$, $s\in (-1,1)$. Then 
\[ 
\fourier{g}(z) = \frac{2\sin z}{z}\quad \text{and}\quad   
f(z) := \fourier{h}(z) = (-2i) \int_0^z \frac{\sin w}{w}\,  dw. 
\]
We claim that
\[ \lim_{\re z \to  \infty, \abs{\im z} < \theta} f(z) = - i\pi
\quad \text{and} \quad
\lim_{\re z \to  -\infty, \abs{\im z} < \theta} f(z) =  i\pi.
\]
This follows from Cauchy's theorem and the well-known fact that
$\int_0^t \sin s / s\, ds \to \pi/2$ as $t \to + \infty$. So in fact
$f \in H^\infty_1(\strip{\theta})$ for each $\theta > 0$. 
\end{exa}

The previous example has an interesting generalisation.
}

In \cite{Haa07b} we considered even functions $g \in \Ell{1}[1,1]$ such 
that, for some $c\in \C$, 
\[ \int_0^1 \abs{ \frac{g(s) - c}{s}}\, ds < \infty;
\]
this reduces the problem to the case that $g$ is constant. Here we
can give a generalisation where we merely assume that $g$ has bounded
variation.

\begin{lemma}\label{ex.l.bv}
Let $g\in \BV[-1,1]$ be an even function and define
$f := \Fourier(\PV-g(s)/s)$. Then  $f \in 
H^\infty_1(\strip{\theta})$, for every $\theta > 0$. 
Moreover, for each $\theta > 0$
there is a constant $c_\theta$ such that
\[ \| f \|_{H^\infty_1({\strip{\theta}})}
\le c_\theta (\Var_{[0,1]}(g) + g(1)) \quad \quad (g \in \BV[0,1],\,\,
f = \Fourier(\PV-g(s)/s)) 
\]
\end{lemma}

\begin{proof}
We first estimate $\abs{zf'(z)}$ and write
\begin{align*}
z f'(z) & =(-i) \fourier{g}(z) = (-2i) \int_0^1 z \cos(sz)g(s)\, ds
= (-2i) \int_0^1 g(s)\, d(\sin(sz))
\\ &  = (-2i)\left( g(1)\sin(z) - \int_0^1 \sin(sz)\, dg(s)\right).
\end{align*}
This yields
\[
\abs{z f'(z)} \le 2 e^\theta (\Var_{[0,1]}(g) + g(1))\quad \quad 
(z\in \strip{\theta}).
\]
To estimate $\abs{f(z)}$ itself, we write
\begin{align*}
f(z) & = (-2i) \int_0^1 \frac{\sin(sz)}{s} g(s)\, ds
= (-2i) \int_0^1 g(s) d\left( \int_0^s \frac{\sin(rz)}{r}\, dr \right)
\\ & = (-2i) \left( g(1)\int_0^1 \frac{\sin(rz)}{r}\, dr
- \int_0^1 \int_0^s \frac{\sin(rz)}{r}\, dr \, dg(s)\right)
\end{align*}
Hence $\abs{f(z)} \le 2c'_\theta (\Var_{[0,1]}(g) + g(1))$, where
$c'_\theta$ is the supremum norm of the function $\int_0^1 \sin(sz)\, ds/s$
on the strip $\strip{\theta}$; this is easily seen to be
finite, cf. \cite[Lemma 3.3]{Haa07b}.
\end{proof}

If $g,f$ are as before, then by Theorem \ref{tra.t.fc}, 
$f(A)$ is a bounded
operator whenever $-iA$ generates a $C_0$-group on a UMD space $X$.
However, we can say (a little) more.

\begin{thm}\label{ex.t.bv}
Let $g \in \BV[-1,1]$ be an even function, and let $f := \Fourier(\PV-g(s)/s)$.
Let $X$ be a UMD space, and let $-iA$ be the generator of a $C_0$-group
$U= (U(s))_{s\in \R}$ on $X$. Then $f(A)$ is bounded, and 
\[ f(A)x \,=\, \PV - \int_{-1}^1 g(s) U(s)x\, \frac{ds}{s}
\,: =\,  \lim_{\epsilon\searrow 0}\,  \int_{\epsilon \le \abs{t}\le 1} 
  g(s) U(s)x\, \frac{ds}{s} \quad\quad (x\in X)
\]
\end{thm}

\begin{proof}
Consider the function $g_\epsilon = g(\car - \car_{(-\epsilon, \epsilon)})$.
Then $g_\epsilon(1) = g(1)$ and $\sup_\epsilon \Var_{[0,1]}(g_\epsilon) < 
\infty$. Let $f_\epsilon := \Fourier(g_\epsilon(s)/s)$. Clearly
$f_\epsilon \to f$ pointwise, and $\sup_\epsilon 
\norm{f_\epsilon}_{H^\infty_1(\strip{\theta})} < \infty$, 
by Lemma \ref{ex.l.bv}. So the statement follows by
Theorem \ref{tra.t.fc} and the Convergence Lemma.
\end{proof}

\vanish{
\begin{rem}
For the special function $g(s) \equiv 1$, one does not need to invoke
the vector-valued Mikhlin result (Lemma \ref{...}). Instead one can argue
directly from the definition of UMD space, using directly the transference
Theorem \ref{...}. The operators involved here are just the 
truncated Hilbert transforms. 
\end{rem}
}

The last class of examples involve functions that are bounded and holomorphic
not only on a strip but on a region
\[ V_{\phi, \theta} := \strip{\theta} \cup \dsector{\phi}
\]
for some $\theta > 0, \phi \in (0, \pi/2)$. (Recall that $\dsector{\phi}$ is
a double sector and $\strip{\theta}$ is a horizontal strip.
The set $V_{\phi, \theta}$  is sometimes called a 
{\em Venturi region} --- inspired by the Venturi tube from fluid dynamics.
We shall need the following fact.

\begin{lemma}\label{ex.l.embed}
Let $0< \theta < \phi < \pi/2$. Then there is a constant $c= c(\theta, \phi)$
such that
\[ \sup\{ \abs{z f'(z)} \suchthat z\in \dsector{\theta}\}
\le c \norm{f}_{H^\infty(\dsector{\phi})} \quad \quad (f\in 
H^\infty(\dsector{\phi})).
\] 
Analogously, for any $0 < \theta < \phi < \pi$
there is a constant $c= c(\theta, \phi)$ such that
\[ \sup\{ \abs{z f'(z)} \suchthat z\in \sector{\theta}\}
\le c \norm{f}_{H^\infty(\sector{\phi})} \quad \quad (f\in 
H^\infty(\sector{\phi})).
\] 
\end{lemma}

\begin{proof}
The proof uses the Cauchy integral integral formula, and is easy.
See \cite[Section 4]{HiePru98} or \cite[Lemma 8.2.6]{HaaFC}.
\end{proof}

Now we can state the theorem, apparently an analogue for unbounded
groups of the Hieber--Pr\"uss Theorem \ref{int.t.hiepru}.

\begin{thm}\label{ex.t.ven}
Let $X$ be a UMD space, and let $-iA$ be the generator
of a strongly continuous group $U= (U(s))_{s\in \R}$ on $X$. 
Let $\theta > \theta(U)$ and $\phi \in (0, \pi/2)$. Then 
$A$ has a bounded $H^\infty(V_{\phi, \theta})$-calculus. More precisely, there
is a constant $c>0$ such that
\[ \norm{f(A + r)} \le c \norm{f}_{\infty} \quad \quad (r\in \R, f\in H^\infty(V_{\phi, \theta}))) 
\]
where $\norm{f}_\infty$ denotes the supremum norm of $f$ on 
$V_{\phi, \theta}$.
\end{thm}

\begin{proof}
By virtue of Theorem \ref{tra.t.fc} 
it suffices to show the continuous inclusion
\[ H^\infty(V_{\phi, \theta}) \subset H^\infty_1(\strip{\theta'})
\] 
for $\theta' \in (\theta(U), \theta)$. But this follows easily
from Lemma \ref{ex.l.embed}, a).
\vanish{Using Cauchy's integral formula
one can show that for each $\phi'\in (0, \phi)$
there is a constant $c$ such that
\[ \abs{z f'(z)} \le c \norm{f}_{\dsector{\phi}, \infty} \quad \quad
(z\in \dsector{\phi'}, f\in H^\infty(\dsector{\phi}))
\]
(see \cite[Lemma 8.2.6]{HaaFC} and compare \cite[Section 4]{HiePru98}). 
Again the Cauchy integral formula
yields a constant $c'$ such that
\[ \abs{z f'(z)} \le c \norm{f}_{\strip{\theta}, \infty}
\quad \quad (z \in \strip{\theta'} \ohne \dsector{\phi'}, f\in H^\infty(\strip{\theta})).
\]
Combining both estimates yields
\[ \abs{z f'(z)} \le \max(c, c')\, \norm{f}_{V_{\phi, \theta}, \infty}
\quad \quad (z\in V_{\phi', \theta'}, f\in H^\infty(V_{\phi, \theta}))
\]
and this is what was to prove.
}
\end{proof}

\vanish{
\begin{rem}\label{ex.r.ven}
Theorem \ref{ex.t.ven} is clearly an analogue for unbounded groups of 
the Hieber--Pr\"uss Theorem \ref{int.t.hiepru}.

of \name{Hieber} and \name{Pr\"uss} from \cite{HiePru98}.
By using the transference principle for bounded groups 
(Theorem \ref{tra.t.bdg}), they showed the following:
If $-iA$ generates a bounded $C_0$-group on a  UMD-space $X$ then $A$ has
a bounded $H^\infty(\dsector{\phi})$-calculus for each $\phi \in (0, \pi/2)$.
See \cite[Theorem 10.7]{KunWei04} for a proof. 
We do not expect a better analogue with respect to $\theta$, 
cf. Remark \ref{tra.r.tra2}.
\end{rem}
}

\section{Sectorial Operators and Cosine Generators}\label{s.cos}

We briefly discuss applications of the previous results to 
sectorial operators and operators that generate cosine functions.

\smallskip

An operator $A$ on a Banach space $X$ is called {\em sectorial}
of angle $\omega \in [0, \pi)$ if 
\[ 
\{ z \in \C \ohne \{0\} \suchthat \omega < \abs{\arg z} \le \pi\} \subset
\resolv(A)
\]
and for every $\omega' \in (\omega, \pi)$
\[ M(A, \omega') := \sup\{ \norm{z R(z,A)} \suchthat
\omega' \le \abs{\arg z} \le \pi\} < \infty.
\]
The minimum  of all $\omega$ such that $A$ is sectorial of angle $\omega$
is denoted by $\omega_{sect}(A)$, and is called the 
{\em sectoriality angle}. 
Basic properties of sectorial operators can be found in 
\cite[Chapter 1]{MarSan} or \cite[Chapter 2]{HaaFC}. 
As in the case of (strong) strip-type operators one has a certain
functional calculus for sectorial operators. (A detailed description 
can be found in \cite{HaaFC}.) If the sectorial operator $A$ is injective, then
$\log A$ is defined, as is  $f(A)$ for each $f\in H^\infty(\sector{\phi})$,
where $\phi \in (\omega, \pi)$. Moreover, $\log(A)$ is a strong strip-type
operator, with $\omega_{sst}(\log(A)) = \omega_{sect}(A)$, 
and there is a composition
rule:
\begin{equation}\label{cos.e.cr}
 f( \log (A)) = (f \nach \log z)(A) \quad \quad 
(\omega \in (\omega_A, \pi), \, f\in H^\infty(\strip{\omega})
\end{equation}
See \cite[Chapter 4]{HaaFC} for these results. The sectorial operator $A$
is said to have {\em bounded imaginary powers} if it is injective
and $-i\log(A)$ generates a $C_0$-group $U$. In this case
$U(s) = A^{-is}$, $s\in \R$. One writes $\theta_A := \theta(U)$ for the
type of this group. By a result of \name{Pr\"uss} and \name{Sohr}
one has  $\omega_{sect}(A) \le \theta_A$, see \cite[Corollary 4.3.4]{HaaFC} or
\cite{Haa03b} for an alternative proof. In the case that $X=H$ is a Hilbert
space then $\theta_A = \omega_{sect}(A)$ 
(a result by \name{McIntosh}, see also 
\cite[Corollary 4.3.5]{HaaFC}).
If $A$ has bounded imaginary powers one writes $A\in \BIP(X)$.

\smallskip

Let $\phi \in (0, \pi)$ and let
\[ H^\infty_{\log}(\sector{\phi})
:= \{ f\in H^\infty(\sector{\phi}) \suchthat
z (\log z) f'(z) \in H^\infty (\sector{\phi})\}
\]
with  the obvious norm. Then we have the following theorem.

\begin{thm}\label{cos.t.mci}
Let $X$ be a UMD space and let $A \in \BIP(X)$ such that $\theta_A < \pi$.
Then the following assertions hold.
\begin{aufzi}
\item $A$ has a bounded $H^\infty_{\log}(\sector{\phi})$-calculus,
for every $\phi \in (\theta_A , \pi)$.
\item If $X = H$ is a Hilbert space, then $A$ has a bounded 
$H^\infty(\sector{\phi})$-calculus, for every $\phi \in (\omega_A, \pi)$.
\end{aufzi}
\end{thm}

Statement b) is due to \name{McIntosh} \cite{McI86}, statement a) is new.
Note that if $f\in H^\infty(\strip{\phi})$ one anyway has that 
$z f'(z)$ is bounded on each smaller sector 
(see the proof of Theorem \ref{ex.t.ven} above). 

\begin{proof}
By virtue of the composition rule (\ref{cos.e.cr}) b) follows from
Corollary \ref{tra.c.bdl} and a) follows from Theorem \ref{tra.t.fc}. 
Note that the mapping
\[ (f \tpfeil f\nach (\log z)) : H^\infty_1(\strip{\phi}) \pfeil
H^\infty_{\log}(\sector{\phi})
\]
is an isometric isomorphism. 
\end{proof}

\medskip

Let us turn to a different application. Note that if $-A$ generates
a bounded $C_0$-semigroup, then $A$ is sectorial of angle $\le \pi/2$.

\begin{thm}\label{cos.t.sgr}
Let $-A$ generate an exponentially stable semigroup $T$ on a UMD space $X$.
If $T$ is a group, then for every $\phi \in (\pi/2, \pi)$ the
$H^\infty(\strip{\phi})$-calculus for $A$ is bounded.
\end{thm}

\begin{proof}
This follows from Theorem \ref{ex.t.ven} by rotating and shifting. 
Furthermore, one needs a certain compatibility of functional calculi
(for sectorial operators and for (rotated, shifted) strip-type operators.
This compatibility is straightforward on the level of elementary calculi
by path deforming, and so holds for extended calculi, cf. 
\cite[Proposition 1.2.7]{HaaFC}.  
\end{proof}

Our techniques allow new proofs of the theorems of \name{Monniaux}
\cite{Mon99} and \name{Dore} and \name{Venni} \cite{DorVen87}; this
is discussed at length in \cite{Haa07b}.

\bigskip

Let us turn to generators of cosine functions. We shall be sketchy
in providing the background, referring to \cite[Section 3.14-3.16]{ABHN}
for the general facts, and to \cite{Haa06fpre} 
for functional calculus matters.

A cosine function on a Banach
space $X$  is a strongly continuous mapping
$\Cos : \R \pfeil \Lin(X)$ that satisfies the identity
\[ \Cos(t +s) + \Cos(t-s) = 2\Cos(t) \Cos(s) \quad \quad (t,s\in \R)
\]
as well as $\Cos(0) = I$.
One can prove from this that a cosine function is exponentially bounded, so
\[ \theta(\Cos) := \inf \big\{ \omega \ge 0 \suchthat \exists M \ge 1 :
\norm{\Cos(t)} \le Me^{\omega \abs{t}}, t \in \R\big\} < \infty. 
\]
The {\em generator} of a cosine function $\Cos$ is defined as the unique
operator $A$ such that 
\begin{equation}\label{cos.e.gen}
 \lambda R(\lambda^2,A) = \int_0^\infty e^{-\lambda t} \Cos(t)\, dt
\quad \quad (\lambda > \theta(\Cos)).
\end{equation}
The cosine function then provides solutions to the
second-order abstract Cauchy problem
\[   u''(t) = Au, \quad \quad u(0)= x, \quad u'(0)=0.
\]
From (\ref{cos.e.gen}) it
 follows that $B:= -A$ is an operator of {\em parabola type} 
$\omega_0:=\theta(\Cos)$, by which 
we mean that $\sigma(B) \subset \cls{\para{\omega_0}}$ and 
for every $\omega > \omega_0$ there exists $M_{\omega}$ such that
\[ \norm{R(\mu, A)} \le \frac{M_\omega}{\sqrt{\abs{\mu}} \big(
\abs{\im \sqrt{\mu}} - \omega\big)}
\quad \quad (\mu \notin \cls{\para{\omega}}).
\]
(Here $\sqrt{\mu}$ denotes any choice of a square root of $\mu$.) As
for strong strip-type or sectorial operators, there is a natural holomorphic
functional calculus associated with such parabola-type operators. The 
procedure is canonical: one considers holomorphic functions $f$ living
on parabolas $\para{\omega}$  with $\omega > \omega_0$. If such a function
has good decay at infinity, one may define
\[ f(B) := \frac{1}{2\pi i} \int_{\rand \para{\omega'}} f(z) R(z,B)\, dz
\]
where $\omega' \in (\omega_0, \omega)$ is arbitrary. This gives
a primary calculus, and by regularisation \cite[Section 1.2]{HaaFC}
one extends this
to a large algebra of meromorphic functions on $\para{\omega}$, including
in particular $H^\infty(\para{\omega})$. If $A= -B$ happens to generate a 
cosine function, one has 
\[ \Cos(t) = \cos(t\sqrt{z})(B) \quad \quad (t\in \R).
\]
Note that since $\cos$ is even, $\cos(t\sqrt{z})$ is a well-defined
bounded holomorphic function on $\para{\omega}$. See \cite{Haa06fpre} 
for proofs and more information.

The idea of reducing the second-order equation to a first-order system
leads to the notion of {\em phase space}. Namely, the operator
matrix 
\[ \calA := \twomat{0}{I}{A}{0}
\]
is the generator of a $C_0$-group $\calU$ on a space of the form
$\calX := X \times V \subset X\times X$, the so-called phase
space. \name{Kisynski} has shown that there is a unique 
subspace of $X\times X$ with
this property, and therefore it was proposed recently that the space
$V$ (which apparently determines the phase space) should be called
the {\em Kisynski space}. It was observed in 
\cite[Appendix]{Haa07b} that $\theta(\calU) = \theta(\Cos)$
and that $V$ and hence $\calX$ is UMD (Hilbert)
if $X$ is UMD  (Hilbert). (This 
was known before, but by an a posteriori argument. See 
\cite{Haa07b} for details on this admittedly cryptic remark.)
As obviously 
\[ \calA^2 = \twomat{A}{0}{0}{A\res{V}} \quad \text{on}\quad
 \DD(A) \times \DD(A\res{V})
\]
is in diagonal form, properties of $A$ can be deduced
from properties of $\calA^2$. In particular, one will get
\[ \twomat{f(-A)}{0}{0}{f(-A)\res{V}} = f(-\calA^2) = f(z^2)(i\calA) 
\]
for sensible holomorphic functions $f$ on the parabola $\para{\omega}$, 
$\omega> \theta(\Cos)= \theta(\calU)$. Writing $\calB := i\calA$
we have that $-i\calB$ generates $\calU$ and we can apply
our results from above. As in the strip case we write
\[ H^\infty_1(\para{\omega}) := \{ f \in H^\infty(\para{\omega}) 
\suchthat z f'(z) \in H^\infty(\para{\omega})\},
\]
and endow it with the canonical norm.

\begin{thm}\label{cos.t.fcpar}
Let $A= -B$ generate a cosine function $\Cos$ on the UMD space $X$, and let
$\omega > \theta(\Cos)$. 
\begin{aufzi}
\item  The operator B has a bounded $H^\infty_1(\para{\omega})$-calculus.

\item If $X= H$ is a Hilbert space, then $B$ has 
a bounded $H^\infty(\para{\omega})$-calculus.
\end{aufzi}
The same statements are true for the operator $B\res{V}$. 
\end{thm}

\begin{proof}
By our remarks above, $f(B) \oplus f(B\res{V})= f(z^2)(\calB)$.  
Now, writing $w = z^2$ we see that 
\[ 2 w f'(w) = z \frac{d}{dz} \big( f(z^2) \big)
\]
so $f\in H^\infty_1(\para{\omega})$ if and only if
$f(z^2) \in H^\infty_1(\strip{\omega})$. Since $-i\calB$ generates
the group $\calU$ on the UMD space $\calX = X \times V$, we can apply
Theorem \ref{tra.t.fc}, a) to conclude that $f(z^2)(\calB)$ is a bounded
operator for all $f\in H^\infty_1(\para{\omega})$, and this proves a).

If $X=H$ is a Hilbert space, then $V$ is also a Hilbert space, and we
may apply Corollary \ref{tra.c.bdl} to prove claim b). 
\end{proof}

\begin{rem}
Part b) of Theorem \ref{cos.t.fcpar} improves the known results
\cite[Section 7.4]{HaaFC} in that we now  have the additional information
that the group on $V \times H$ has the same growth type as the original
cosine function, and so no shifting is needed any more.
\end{rem}

\medskip

If $B$ is an operator of parabola-type $\omega_0$, then for large
$\lambda > 0$ the operator $\lambda + B$
will be sectorial. In fact, simple geometry yields that
for $\theta \in (0, \pi/2]$ and $\omega \ge 0$ 
\begin{equation}\label{cos.e.inc}
 \left(\mfrac{\omega}{\sin \theta} \right)^2 + \para{\omega} 
\subset S_{\theta},
\end{equation}
and some further computation shows that the operator
$B_\theta := B + (\omega_0/ \sin\theta)^2$ is sectorial of angle 
$\theta$, see \cite[Proposition 7.6]{Haa06fpre}. 
Furthermore, if $-B$ generates a cosine function then so does $-B_\theta$,
by perturbation theory \cite[Corollary 3.14.10]{ABHN}.

\begin{thm}\label{cos.t.fcsec}
Let $A= -B$ be the generator of a cosine function $\Cos$
on the UMD space $X$.
Let $\theta \in (0,\pi/2]$ and set $B_\theta := -A + (\omega_0/\sin\theta)^2$,
where $\omega_0 := \theta(\Cos)$ is the exponential growth type
of the cosine function. Then 
the operator $B_\theta$ has bounded  $H^\infty(S_\phi)$-calculus
for every $\phi \in (\theta, \pi)$.
\end{thm}

\begin{proof}
Choose $\theta' \in (\theta, \phi)$ and define
$\omega:= \omega_0 \sin\theta'/\sin\theta > \omega_0$ and 
$\lambda := (\omega/\sin\theta')^2 = (\omega_0/\sin\theta)^2$.
Let $f\in H^\infty(S_\phi)$ and define 
$g (z) := f( \lambda +z)$.
Then obviously 
$f(B_\theta) = f( (\omega_0/\sin\theta)^2 + B) =
f( \lambda + B) = g(B)$. Now $g \in H^\infty_1(\para{\omega})$.
To see this note first that
by (\ref{cos.e.inc})
\[
 \lambda + \para{\omega} = 
\left(\mfrac{\omega}{\sin \theta'} \right)^2 + \para{\omega} 
\subset S_{\theta'} \subset S_\phi
\]
so $g$ is bounded on $\para{\omega}$ by $\norm{f}_{H^\infty(S_\phi)}$.
Moreover,
\[
\abs{zg'(z)} =  \abs{\frac{z}{z + \lambda}}
\abs{(z + \lambda) f'(z + \lambda)}
\le  c \abs{\frac{z}{z + \lambda}} \norm{f}_{H^\infty(S_{\phi})}
 \quad \quad (z\in \para{\omega})
\]
for some constant $c$, by Lemma \ref{ex.l.embed}; since
$-\lambda \notin \cls{\para{\omega}}$, the first factor
is bounded on $\para{\omega}$. 
We may now apply Theorem \ref{cos.t.fcpar} to conclude that $f(B_\theta)
= g(B)$ is bounded.
\end{proof}

\begin{cor}\label{cos.c.fcsec}
Let $A$ be the generator of a cosine function on the UMD space $X$.
If $B:= -A$ is sectorial and invertible, then $B$ has
bounded $H^\infty(S_\phi)$-calculus for every 
$\phi \in (\omega_{sect}(B), \pi)$.
\end{cor}

\begin{proof}
Fix $\phi \in (\omega_{sect}(B), \pi)$. 
Since $B$ is assumed to be invertible,
standard perturbation theory \cite[Corollary 5.5.5]{HaaFC} shows that it
suffices to prove that $\lambda + B$ has
bounded $H^\infty(S_\phi)$-calculus, for some $\lambda>0$.
By Theorem \ref{cos.t.fcpar}, $\lambda := (\omega_0/\sin\theta)^2$
will do, where $\omega_0$ is the exponential growth type of the
cosine function generated by $A= -B$ and $\theta\in (0, \phi)$ is arbitrary.
\end{proof}

\begin{rem}
In analogy with groups one would expect much stronger results
for {\em bounded} cosine functions. And indeed, in \cite{Haa08b} 
a transference principle for bounded cosine functions was 
established, and it was used to show that every bounded
cosine function on a UMD space has a uniformly bounded
square root reduction group. This had been an open problem
for quite some time.
\end{rem}

\bibliographystyle{plain}

\def\cprime{$'$} \def\cprime{$'$} \def\cprime{$'$} \def\cprime{$'$}


\end{document}